\documentclass[11pt,a4paper]{article}
\usepackage[margin=1in]{geometry}
\usepackage[T1]{fontenc}
\usepackage{lmodern}
\usepackage{amsmath,amssymb,amsthm,mathtools}
\usepackage{microtype}
\usepackage{enumitem}
\usepackage{xcolor}
\usepackage{hyperref}
\hypersetup{colorlinks=true,linkcolor=blue!45!black,
  citecolor=blue!45!black,urlcolor=blue!45!black,
  pdftitle={Spectral Extremal Problems for Graphs with Prescribed Integrity}}
\setlist{itemsep=2pt,topsep=4pt}
\newtheorem{theorem}{Theorem}[section]
\newtheorem{lemma}[theorem]{Lemma}
\newtheorem{proposition}[theorem]{Proposition}
\newtheorem{corollary}[theorem]{Corollary}
\theoremstyle{remark}
\newtheorem{remark}[theorem]{Remark}
\newcommand{\I}{\mathrm{I}}
\newcommand{\Spec}{\operatorname{Spec}}
\newcommand{\mult}{\operatorname{mult}}
\newcommand{\ones}{\mathbf{1}}
\newcommand{\R}{\mathbb{R}}

\title{Signless Laplacian Spectral Extremal Problems\\
for Graphs with Prescribed Integrity and Bounded Connectivity}
\author{Yamini V$^*$ \quad and \quad Joyentanuj Das\footnote{Department of Mathematics, College of Engineering and Technology, SRM Institute of Science and Technology, Kattankulathur, Chennai, 603203, India. \indent  Email: joyentanuj@gmail.com,  joyentad@srmist.edu.in}}
\date{}

\begin{document}
\maketitle

\begin{abstract}
Let $\I(G)=\min_{S\subseteq V(G)}\{|S|+\psi(G-S)\}$ denote the
vertex integrity of a graph $G$, where $\psi(G-S)$ is the order of
a largest component of $G-S$. We determine the unique graph maximizing
the signless Laplacian spectral radius among connected graphs of order
$n$ with integrity at most $t$, and also under the simultaneous
constraint that the vertex-connectivity is at most $k$. In the
unrestricted case the extremal graph is
$F_{n,t}=K_{t-1}\vee\overline K_{n-t+1}$. If
$1\le k<t-1\le n-2$, the unique extremal graph under the additional
connectivity constraint is
\[
H_{n,k,t}=K_k\vee\bigl(K_1\cup
(K_{t-k-1}\vee\overline K_{n-t})\bigr).
\]
Both extremal graphs have integrity exactly $t$, so the results also
solve the corresponding problems with integrity prescribed exactly.
The main argument combines integrity and cut certificates with
Perron-vector edge relocations that preserve both constraints.
We give quotient-matrix descriptions of the relevant signless
Laplacian spectral radii, determine the Laplacian spectra of the
extremal constructions, and derive fixed-parameter asymptotic
expansions for their adjacency spectral radii.
\end{abstract}

\noindent\textbf{Keywords:} vertex integrity; signless Laplacian;
spectral extremal graph; vertex-connectivity; Laplacian spectrum.

\noindent\textbf{2020 Mathematics Subject Classification:}
05C50, 05C40.

\section{Introduction and preliminaries}\label{sec:introduction}

The integrity of a graph measures the combined cost of deleting
vertices and controlling the order of the largest remaining component.
It was introduced in the study of graph vulnerability by
Barefoot, Entringer, and Swart \cite{Barefoot1987}.
For a finite simple graph $G$, put
\begin{equation}\label{eq:integrity}
 \I(G)=\min_{S\subseteq V(G)}
 \bigl(|S|+\psi(G-S)\bigr),
\end{equation}
where $\psi(X)$ is the order of a largest connected component of $X$,
and $\psi(\varnothing)=0$. A set attaining the minimum in
\eqref{eq:integrity} is an \emph{integrity set}.
Throughout, graphs are undirected and have no loops or multiple edges.

Determining extremal spectral radii under structural constraints
is a well-established direction in spectral graph theory. Brualdi
and Solheid \cite{BrualdiSolheid1986} gave an influential early
framework for extremal spectral-radius questions. For the signless
Laplacian, the standard reference is the work of Cvetkovi\'c,
Rowlinson, and Simi\'c \cite{Cvetkovic2007}; see also the monograph
of Brouwer and Haemers \cite{BrouwerHaemers2012}. Ye, Fan, and Wang
\cite{YeFanWang2010} determined extremal adjacency and signless
Laplacian spectral radii under fixed connectivity constraints.
These results provide a natural background for studying spectral
extremal problems in which a global vulnerability parameter and a
connectivity restriction are imposed together.

Several related extremal spectral problems have also been
investigated by the authors and their collaborators. In
particular, Das and Mahato \cite{DasMahato2024} studied the
$A_{\alpha}$-spectral radius of unicyclic and bicyclic graphs
under prescribed girth and pendant-vertex constraints. This
work is particularly relevant to the present investigation
because $A_{1/2}(G)=\frac{1}{2}Q(G)$, and hence the signless
Laplacian spectral radius is twice the $A_{1/2}$-spectral
radius. In a different direction, Das \cite{DasCliqueTrees2025}
determined the extremal spectral radius of clique trees with a
prescribed zero forcing number. Das and Mohanty
\cite{DasMohantyBlock2024} studied the maximization of the
spectral radius of block graphs with a given dissociation
number. These works demonstrate the use of structural and
spectral techniques in extremal problems involving prescribed
graph parameters. The present work develops this line of
investigation further by considering the signless Laplacian
spectral radius under simultaneous restrictions on vertex
integrity and vertex-connectivity.

The closest recent work to the present problem is that of Dai, Shen,
Song, and Liu \cite{DaiShenSongLiu2027}, who studied the maximum
spectral radius of graphs with prescribed integrity and bounded
connectivity. Their paper makes the interaction between integrity
and vertex-connectivity explicit in an adjacency-spectral extremal
setting. The purpose here is different: we investigate the
signless Laplacian spectral radius, whose degree contribution changes
both the extremal comparisons and the appropriate Perron-vector
relocation argument. Thus the present results should be read as a
signless-Laplacian counterpart, not as a replacement for the
adjacency-spectral problem treated in \cite{DaiShenSongLiu2027}.
The distinction is important because an extremal theorem for the
adjacency matrix does not automatically imply the corresponding
statement for the signless Laplacian.

Write $A(G)$ and $D(G)$ for the adjacency and diagonal degree matrices,
and put
\[
 L(G)=D(G)-A(G),\qquad Q(G)=D(G)+A(G).
\]
Their largest eigenvalues are denoted by $\rho(G)$, $\mu_1(G)$,
and $q(G)$, respectively. The complement of $G$ is $\overline G$.
The join of vertex-disjoint graphs $X,Y$, obtained by adding every
edge between their vertex sets, is denoted by $X\vee Y$.
The complete and edgeless graphs of order $a$ are $K_a$ and
$\overline K_a$. Empty parts are allowed when explicitly indicated.
The vertex-connectivity of a connected noncomplete graph is the
minimum number of vertices whose deletion disconnects it; as usual,
$\kappa(K_n)=n-1$.

For $2\le t\le n$, define
\begin{equation}\label{eq:F}
 F_{n,t}=K_{t-1}\vee\overline K_{n-t+1}.
\end{equation}
For $1\le k<t-1\le n-2$, define
\begin{equation}\label{eq:H}
 H_{n,k,t}=K_k\vee
 \bigl(K_1\cup(K_{t-k-1}\vee\overline K_{n-t})\bigr).
\end{equation}
In the second construction the outer $K_k$ is universal. The
distinguished vertex in $K_1$ is adjacent precisely to that clique.
Equivalently, $H_{n,k,t}$ is obtained from $F_{n-1,t}$ by adding one
vertex adjacent to $k$ vertices of its $(t-1)$-clique.

Our principal results are the inequalities
\[
 q(G)\le q(F_{n,t})\quad\text{if }\I(G)\le t,
\]
and
\[
 q(G)\le q(H_{n,k,t})
 \quad\text{if }\I(G)\le t\text{ and }\kappa(G)\le k,
\]
for connected graphs of order $n$, with equality precisely for the
displayed graphs. The second assertion is the principal result of the paper. Its
proof simultaneously controls an integrity certificate and a vertex cut. Working first with
$\I(G)\le t$ allows certificate-preserving transformations even when
they lower the integrity. The equality examples have integrity
exactly $t$.

The Laplacian statements later in the paper are applications of the classical complement and join identities, which we prove in detail and attribute to the standard literature \cite{BrouwerHaemers2012,Mohar1991}. The adjacency calculations give exact formulas and asymptotic expansions for the explicit constructions \eqref{eq:F}--\eqref{eq:H}; they are not asserted to solve the adjacency-spectral extremal problem considered in
\cite{DaiShenSongLiu2027}.

In the next lemma, we prove some of the observations that will be used later.

\begin{lemma}\label{lem:integrity}
The following statements hold.
\begin{enumerate}[label=(\roman*)]
\item If $G,G'$ have the same vertex set and $E(G)\subseteq E(G')$,
then $\I(G)\le\I(G')$.
\item For every graph $X$ and integer $a\ge1$,
$\I(K_a\vee X)=a+\I(X)$.
\item If $X$ is nonempty, then $\I(X\cup K_1)=\I(X)$.
\end{enumerate}
\end{lemma}
\begin{proof}
For every vertex set $S$, each component of $G-S$ is contained
in a component of $G'-S$. Hence
$\psi(G-S)\le\psi(G'-S)$, and taking minima proves (i).

For (ii), delete the complete part together with an integrity set of
$X$ to obtain the upper bound. Conversely, let $S$ be any vertex set
of $K_a\vee X$. If some vertex of the complete part remains, the
remaining graph is connected, so the value of the expression in
\eqref{eq:integrity} is $a+|V(X)|\ge a+\I(X)$. If the entire
complete part is deleted, write $S=V(K_a)\cup S_X$. The same
expression is
$a+|S_X|+\psi(X-S_X)\ge a+\I(X)$. This proves (ii), including
the case $X=\varnothing$.

For (iii), let $b$ be the added isolated vertex and put
$S_X=S\cap V(X)$. If $b\notin S$, then
\[
 |S|+\psi((X\cup K_1)-S)
 =|S_X|+\max\{\psi(X-S_X),1\}\ge\I(X).
\]
If $b\in S$, the corresponding expression is
$1+|S_X|+\psi(X-S_X)\ge\I(X)$. For the reverse inequality,
choose an integrity set $S_X$ leaving at least one vertex of $X$.
Such a choice is always possible: if deleting all vertices is
optimal, deleting all but one has the same value $|V(X)|$.
Leaving $b$ undeleted now gives exactly
$|S_X|+\psi(X-S_X)=\I(X)$.
\end{proof}

\begin{proposition}\label{prop:feasibility}
For the parameter ranges in \eqref{eq:F} and \eqref{eq:H},
\[
 \I(F_{n,t})=t,\qquad
 \I(H_{n,k,t})=t,\qquad
 \kappa(H_{n,k,t})=k.
\]
\end{proposition}
\begin{proof}
A nonempty edgeless graph has integrity one, so
Lemma~\ref{lem:integrity}(ii) gives
$\I(F_{n,t})=(t-1)+1=t$. Also
\[
 H_{n,k,t}=K_k\vee(K_1\cup F_{n-k-1,t-k}).
\]
The inner split graph is nonempty and has integrity $t-k$.
Parts (ii) and (iii) of the lemma give $\I(H_{n,k,t})=t$.
Deleting the universal $k$-clique disconnects the graph. Deleting
fewer than $k$ vertices leaves a universal vertex, and hence leaves
a connected graph. Therefore its vertex-connectivity is exactly $k$.
\end{proof}

For every vector $x\in\R^{V(G)}$,
\begin{equation}\label{eq:Qform}
 x^{\mathsf T}Q(G)x=\sum_{uv\in E(G)}(x_u+x_v)^2.
\end{equation}
If $G$ is connected, $Q(G)$ is irreducible and nonnegative.
Perron--Frobenius theory therefore provides a strictly positive
eigenvector for $q(G)$; the eigenspace is one-dimensional.
We use the Rayleigh principle and these standard facts as in
\cite{BrouwerHaemers2012,Cvetkovic2007}.
In particular, adding a missing edge to a connected graph strictly
increases $q$, by evaluating \eqref{eq:Qform} at a positive unit
Perron vector of the original graph.

A partition is \emph{equitable} for a matrix when each of its blocks
has a constant row sum. Its quotient matrix records those row sums.
The space of vectors constant on each cell is invariant, and the
quotient represents the action on that space. If the matrix and
quotient are nonnegative and irreducible, a positive quotient
eigenvector lifts to a positive eigenvector of the full matrix.
Their Perron eigenvalues therefore coincide. We will also use the
following immediate consequence of uniqueness: a positive unit
Perron vector is constant on each orbit of graph automorphisms.

The paper is organized as follows. In Section~2, we investigate the
unrestricted signless Laplacian extremal problem and establish that
$F_{n,t}=K_{t-1}\vee\overline K_{n-t+1}$ is the unique graph attaining
the maximum signless Laplacian spectral radius among connected graphs
of order $n$ with vertex integrity at most $t$. In Section~3, we
consider the simultaneous restrictions on vertex integrity and
vertex-connectivity. This section contains our principal result,
which characterizes the unique extremal graph
$H_{n,k,t}=K_k\vee(K_1\cup(K_{t-k-1}\vee\overline K_{n-t}))$
under these constraints. In Section~4, we determine the Laplacian
spectra of the extremal constructions and establish the corresponding
spectral extremality results. Section~5 is devoted to the adjacency
spectral properties of these graphs, where we derive exact expressions
for their adjacency spectral radii and obtain asymptotic expansions
when the integrity parameters are fixed. Finally, Section~6 concludes
the paper with a summary of the main findings and some directions for
further research.

\section{The unrestricted signless Laplacian problem}\label{sec:unrestricted}

\begin{proposition}\label{prop:qF}
For $2\le t\le n$,
\begin{equation}\label{eq:qF}
 q(F_{n,t})=
 \frac{n+2t-4+
 \sqrt{(n+2t-4)^2-8(t-1)(t-2)}}{2}.
\end{equation}
\end{proposition}
\begin{proof}
Put $r=t-1$ and $p=n-t+1$. The vertices of the clique have
degree $n-1$, and the independent vertices have degree $r$.
Thus the equitable quotient of $Q(F_{n,t})$ is
\[
 \begin{pmatrix}n+r-2&p\\r&r\end{pmatrix}.
\]
It is nonnegative and irreducible, and its characteristic polynomial
is
$x^2-(n+2t-4)x+2(t-1)(t-2)$.
Its larger root is the Perron eigenvalue, giving \eqref{eq:qF}.
\end{proof}

\begin{theorem}\label{thm:unrestricted-q}
Let $G$ be a connected graph of order $n$ and let $2\le t\le n$.
If $\I(G)\le t$, then
\[
 q(G)\le q(F_{n,t}),
\]
with equality if and only if $G\cong F_{n,t}$.
In particular, the same graph is the unique maximizer among connected
graphs of order $n$ and integrity exactly $t$.
\end{theorem}
\begin{proof}
If $t=n$, every graph is a spanning subgraph of
$K_n=F_{n,n}$, and strict edge monotonicity gives the conclusion.
Assume $t<n$. Choose a set $S$ certifying $\I(G)\le t$, and put
$s=|S|$ and $m=t-s$. Since $G$ is connected and $t<n$,
\[
 1\le s\le t-1,\qquad m\ge1.
\]
Let the components of $G-S$ have orders $c_1,\ldots,c_\ell$.
Then $1\le c_i\le m$ and $\sum_i c_i=n-s$.
Completing $S$ and each component, and adding all edges between
$S$ and its complement, gives
\[
 G^*=K_s\vee(K_{c_1}\cup\cdots\cup K_{c_\ell}).
\]
The same set $S$ certifies $\I(G^*)\le t$.
Moreover $q(G)\le q(G^*)$, with strict inequality unless
no edge was added.

Write $\lambda=q(G^*)$. A positive Perron vector has a constant
coordinate $x_0$ on $S$ and a constant coordinate $x_i$ on the
$i$th component. The eigenvector equations are
\begin{align*}
 \lambda x_0&=(n+s-2)x_0+\sum_i c_i x_i,\\
 \lambda x_i&=s x_0+(s+2c_i-2)x_i.
\end{align*}
In particular, each denominator in the identity
\begin{equation}\label{eq:q-completed}
 \lambda=n+s-2+
 s\sum_i\frac{c_i}{\lambda-s-2c_i+2}
\end{equation}
is positive.

Put $Q=q(F_{n,t})$. We first record
\begin{equation}\label{eq:Qpositive}
 Q>2t-2.
\end{equation}
Indeed, $F_{n,t}$ contains a $K_t$. Extend a positive unit Perron
vector of that clique by zeros on the other vertices.
Its Rayleigh quotient in $F_{n,t}$ exceeds $2t-2$, because
$n>t$ and there is an edge from the clique to an additional
independent vertex. This proves \eqref{eq:Qpositive}.
Consequently $Q+s-2t+2>0$.

For $1\le s\le t-1$, define
\[
 P_s(x)=(x-n-s+2)(x+s-2t+2)-s(n-s).
\]
Expanding gives
\[
 P_s(x)=x^2-(n+2t-4)x
       +2\bigl((t-1)(n-2)-s(n-t)\bigr).
\]
The quotient polynomial in Proposition~\ref{prop:qF} is
$P_{t-1}$, so $P_{t-1}(Q)=0$ and
\begin{equation}\label{eq:Ps}
 P_s(Q)=2(t-1-s)(n-t)\ge0.
\end{equation}
Since $Q+s-2t+2>0$, this is equivalent to
\begin{equation}\label{eq:comparisonQ}
 n+s-2+\frac{s(n-s)}{Q+s-2t+2}\le Q.
\end{equation}

Suppose that $\lambda>Q$. Then
$\lambda+s-2t+2>0$, and $c_i\le t-s$ gives, from
\eqref{eq:q-completed},
\[
 \lambda\le n+s-2+
 \frac{s(n-s)}{\lambda+s-2t+2}
 < n+s-2+\frac{s(n-s)}{Q+s-2t+2}\le Q.
\]
The middle inequality uses $s(n-s)>0$. This contradiction proves
$q(G)\le\lambda\le Q$. This argument uses the denominator comparison
only when $\lambda\ge Q$, where its positivity has been established.

If $q(G)=Q$, then $\lambda=Q$. When $s<t-1$,
\eqref{eq:Ps} and \eqref{eq:comparisonQ} are strict, and the same
denominator comparison at $\lambda=Q$ is impossible.
Therefore $s=t-1$, so $m=1$ and $G^*=F_{n,t}$.
Strict edge monotonicity forces $G=G^*$.
Conversely, Proposition~\ref{prop:feasibility} proves feasibility
of $F_{n,t}$, and equality is immediate.
\end{proof}

\section{A simultaneous integrity and connectivity bound}\label{sec:bounded}

The difficulty in imposing a connectivity bound is that completing
edges relative to an integrity set may destroy a small vertex cut.
We avoid this by exchanging the two certificates before completing
the graph. The edge-relocation lemma below includes the case of equal
Perron coordinates, which is needed for uniqueness.

\begin{lemma}[Strict signless-Laplacian edge relocation]\label{lem:q-relocation}
Let $G$ be connected, let $x$ be a positive unit Perron vector of $Q(G)$,
and let $a,b$ be distinct vertices with $x_a\ge x_b$. Let
$W\subseteq N_G(b)\setminus(N_G(a)\cup\{a\})$ be nonempty. Form $G'$
by replacing $bz$ with $az$ for every $z\in W$. Then $q(G')>q(G)$.
\end{lemma}
\begin{proof}
The signless Laplacian quadratic form gives
\[
x^{\mathsf T}(Q(G')-Q(G))x
=(x_a-x_b)\sum_{z\in W}(x_a+x_b+2x_z)\ge0.
\]
If $x_a>x_b$, the inequality is strict. If $x_a=x_b$ and the two
spectral radii were equal, equality in the Rayleigh principle would imply
$Q(G')x=q(G)x$. But at the coordinate $b$,
\[
(Q(G')x)_b=(Q(G)x)_b-|W|x_b-\sum_{z\in W}x_z<q(G)x_b,
\]
a contradiction. This proves strictness also when the two Perron
coordinates are equal.
\end{proof}

\begin{remark}\label{rem:q-multirelocation}
The same argument applies to several donors $b_1,\ldots,b_\ell$ and one
receiver $a$, provided that the moved edges and the target edges are each pairwise
distinct, every target edge is absent, $x_a\ge x_{b_i}$ for every donor, and at least one donor loses
an edge and receives no new edge. The quadratic-form differences add.
If their sum vanishes, the eigenvector equation fails at any such donor.
We use this observation only for donors in pairwise disjoint cliques,
with all target edges incident to a vertex in a different clique.
\end{remark}

\begin{theorem}\label{thm:bounded-q-solved}
Let $1\le k<t-1\le n-2$. Among all connected graphs $G$ of order $n$
with $\I(G)\le t$ and $\kappa(G)\le k$, the unique graph maximizing
$q(G)$, up to isomorphism, is
\[
H_{n,k,t}=K_k\vee\bigl(K_1\cup(K_{t-k-1}\vee\overline K_{n-t})\bigr).
\]
Consequently, the same conclusion holds under the exact constraint
$\I(G)=t$.
\end{theorem}
\begin{proof}
The candidate has $\I(H_{n,k,t})=t$ and $\kappa(H_{n,k,t})=k$, so the
class is nonempty. Choose a graph $G$ in this finite class with maximum
signless Laplacian spectral radius. Write $x$ for a positive Perron
vector of $Q(G)$. An \emph{integrity certificate} will mean a set
$S\subseteq V(G)$ such that $|S|+\psi(G-S)\le t$; equivalently,
$|S|\le t$ and every component of $G-S$ has order at most
$m=t-|S|$. Since $t<n$ and $G$ is connected, every such set satisfies
$1\le |S|\le t-1$, and $G-S$ has at least two components. A
\emph{cut certificate} will mean a vertex cut $T$ with
$1\le |T|\le k$.

For any fixed pair of certificates $S,T$, an absent edge can be added
whenever it is allowed by both certificates: its endpoints either meet
$S$ or belong to one component of $G-S$, and they either meet $T$ or
belong to one component of $G-T$. Such an addition leaves both
certificates valid and strictly increases $q$. By extremality, every
edge allowed by both certificates is therefore already present. We
refer to this fact as saturation.

Among all pairs of certificates for $G$, choose $S,T$ with
$|S\cap T|$ maximum. Suppose that neither contains the other. Choose
$u\in T\setminus S$ and $v\in S\setminus T$. Let $C$ be the component
of $G-S$ containing $u$, and let $D$ be the component of $G-T$
containing $v$. Saturation gives $uv\in E(G)$, since this edge meets
both $S$ and $T$.

If $x_v\ge x_u$, put
\[
W=N_G(u)\setminus(T\cup V(D)).
\]
For each $z\in W$, the edge $vz$ is absent, because $v$ and $z$ lie in
different components of $G-T$. Replace all $uz$, $z\in W$, by $vz$.
All new edges meet $S$, so $S$ remains an integrity certificate. Set
$T'=(T\setminus\{u\})\cup\{v\}$. In the graph obtained after deleting
$T'$, the vertices of $(V(D)\setminus\{v\})\cup\{u\}$ have no edges
to the other original components of $G-T$. This set is nonempty,
as are those other components, so $T'$ is a cut certificate. The
modified graph is connected: each removed edge $uz$ is replaced by
the path $u,v,z$, and the edge $uv$ is retained. If $W$ is nonempty,
Lemma~\ref{lem:q-relocation} contradicts extremality. If $W$ is empty,
the graph has not changed and $|S\cap T'|=|S\cap T|+1$, contradicting
the choice of the certificates.

If $x_u\ge x_v$, the symmetric operation uses
\[
W=N_G(v)\setminus(S\cup V(C)).
\]
Now every target edge $uz$ is absent because $u,z$ lie in distinct
components of $G-S$. Replace $vz$ by $uz$ for $z\in W$, retain $T$,
and set $S'=(S\setminus\{v\})\cup\{u\}$. The only possible new
component of the residual graph containing $v$ lies in
$(V(C)\setminus\{u\})\cup\{v\}$, whose order is $|C|\le m$;
all other residual components are contained in the old ones. Thus
$S'$ is an integrity certificate. All new edges meet $T$, so the
original cut remains valid. Connectivity follows from the retained
edge $uv$, just as above. A nonempty $W$ contradicts extremality, and
an empty $W$ increases $|S\cap T|$ by one. These contradictions prove
that the certificates are nested.

If $S\subseteq T$, then $|S|\le k$ and $S$ itself is a cut certificate,
since $n-|S|>t-|S|$. Replace $T$ by $S$. We have therefore arranged
\[
T\subseteq S,\qquad h=|T|\le k,\qquad s=|S|,\qquad m=t-s.
\]
Saturation implies that $T$ is a universal clique. In each component
$D_i$ of $G-T$, the set $S\cap V(D_i)$ is a clique adjacent to every
other vertex of $D_i$, and the components of $G-S$ contained in $D_i$
are cliques of order at most $m$.

First suppose that $s>h$. There exist vertices
$u\in S\setminus T$ and $v\in V(G)\setminus S$ in different components
of $G-T$. Indeed, if every such pair belonged to one component, fix
one vertex from each nonempty set; all vertices outside $T$ would then
belong to that same component, contradicting that $T$ is a cut.
Let $a$ be whichever of $u,v$ has the larger Perron coordinate, and
let $b$ be the other. Transfer all edges $bz$ with $z\notin T$ to
$az$. The target edges are absent, because $a,b$ belong to distinct
components of $G-T$. The universal clique $T$ remains unchanged, so
the modified graph is connected. The vertex $b$ has precisely the
vertices of $T$ as neighbours, and hence $T$ is still a cut.

If $a\in S$, retain $S$; the new edges meet $S$ and $b\notin S$. If
$b\in S$, replace it in the certificate by $a$. In the latter case,
the old component of $G-S$ containing $a$ loses $a$, the new residual
vertex $b$ is isolated, and all other residual components only lose
edges. Since $m\ge1$, the certificate remains valid. Thus the new
certificate, which we again call $S$, has size $s$, contains $T$, and
excludes $b$. A nonempty transfer strictly increases $q$, which is
impossible at a maximizer. Consequently the transfer was empty, the
graph remains $G$, and we have obtained the stated singleton cut
component merely by changing the certificate if necessary.

Complete all edges on $V(G)\setminus\{b\}$ allowed by this integrity
certificate, and retain exactly the edges from $b$ to $T$. This
preserves the certificate and the cut isolating $b$. Extremality and
strict edge monotonicity imply that these edges were already present.
Consequently $G$ has the form
\begin{equation}\label{eq:singleton-normal-form}
K_h\vee\left(K_1\cup\left(K_{s-h}\vee
\bigcup_{i=1}^{\ell}K_{c_i}\right)\right),
\qquad 1\le c_i\le m,\qquad \sum_i c_i=n-s-1.
\end{equation}

It remains in this step to treat $s=h$. In this case saturation gives
\[
G=K_h\vee\bigcup_{i=1}^{\ell}K_{c_i},\qquad
1\le c_i\le m,\qquad \ell\ge2.
\]
Here $m=t-h\ge t-k\ge2$. Order the cliques so that
$c_1\ge c_i$ for every $i$, choose $a\in K_{c_1}$, and choose
$b\in K_{c_2}$. By symmetry the Perron coordinate on a clique of
order $c$ is
\[
\frac{\sum_{z\in T}x_z}{q(G)-h-2c+2}.
\]
Every denominator is positive by the corresponding eigenvector
equation, and this expression increases with $c$. Hence $x_a$
is at least every Perron coordinate outside $T$.

Transfer the edges from $b$ to $K_{c_2}\setminus\{b\}$ to $a$.
For every $i\ge3$ with $c_i=m$, select a vertex $v_i\in K_{c_i}$ and
transfer its $m-1$ internal edges to $a$. All target edges are absent,
and the receiver dominates the donors in the original Perron vector.
Next add all missing edges from $a$ to vertices outside $T$, except
for $b$. The resulting graph has integrity certificate
$S'=T\cup\{a\}$ with capacity $m-1$: the first clique loses $a$,
the second clique loses $b$, every other clique of size $m$ loses its
selected vertex, and each selected vertex is isolated in the residual
graph. The vertex $b$ remains adjacent only to $T$, so the cut is
preserved. The intermediate transfer graph has the same certificates,
and all these graphs are connected because $T$ stays universal.
Thus the construction stays in the feasible class throughout. If
any edge is moved, the simultaneous version of
Lemma~\ref{lem:q-relocation} gives a strict increase, a contradiction.
No edge addition can be strict at the maximizer either. We have therefore
obtained \eqref{eq:singleton-normal-form}, now with $s=h+1$ and
capacity $t-s=m-1$, without changing $G$.

Suppose that $G$ has the form \eqref{eq:singleton-normal-form} and
$m=t-s\ge2$. Since
\[
\sum_i c_i=n-s-1=m+(n-t)-1\ge m,
\]
there is at least one residual clique. Choose $a$ in a largest one,
say $K_{c_1}$. The Perron coordinate on a residual clique of order
$c$ is
\[
\frac{\sum_{z\in S}x_z}{q(G)-s-2c+2},
\]
so $x_a$ is at least the coordinate on any other residual clique.
For each $i\ge2$ with $c_i=m$, choose $v_i\in K_{c_i}$ and transfer
every edge from $v_i$ to $K_{c_i}\setminus\{v_i\}$ to $a$.
Then add all missing edges from $a$ to
$V(G)\setminus(T\cup\{a,b\})$. Equivalently, make $a$ universal in
the core $G-b$, while leaving $ab$ absent.

The new certificate $S'=S\cup\{a\}$ has capacity $m-1$.
The first residual clique loses $a$, each other clique of order $m$
loses its selected vertex, the selected vertices become residual
singletons, and all remaining cliques already have order at most
$m-1$. If no clique has order $m$, the same conclusion holds without
any transfers. The vertex $b$ still has neighbourhood $T$, and $T$
is universal. Thus connectedness and the small cut are preserved.
Every nonempty transfer or edge addition strictly increases $q$;
at the maximizer the graph is unchanged. The new normal form has
$s$ replaced by $s+1$ and $m$ by $m-1$.

Iterate this step until $m=1$. Then $s=t-1$, and
\eqref{eq:singleton-normal-form} is exactly $H_{n,h,t}$. Finally,
if $h<k$, choose $k-h$ vertices from the inner clique of order
$t-h-1$ and join them to $b$. This produces $H_{n,k,t}$ by strictly
adding edges; it is feasible because its integrity is $t$ and its
connectivity is $k$. This contradicts extremality unless $h=k$.
We conclude that every maximizer is isomorphic to $H_{n,k,t}$.
Since the candidate belongs to the class, the asserted maximum and
the unique equality case follow. The exact-integrity class is a
subclass containing the same candidate, proving the final assertion.
\end{proof}

\begin{proposition}\label{prop:qH}
Let $1\le k<t-1\le n-2$. Then $q(H_{n,k,t})$ is the largest
eigenvalue of
\begin{equation}\label{eq:N}
 N_{n,k,t}=
 \begin{pmatrix}
 n+k-2&1&t-k-1&n-t\\
 k&k&0&0\\
 k&0&n+t-k-4&n-t\\
 k&0&t-k-1&t-1
 \end{pmatrix}.
\end{equation}
\end{proposition}
\begin{proof}
Partition the graph into the universal clique $Z$, the distinguished
singleton $\{b\}$, the inner clique $A$, and the independent set $P$.
Their sizes are $k,1,a=t-k-1,p=n-t$, and their degrees are
$n-1,k,n-2,t-1$, respectively.
The row sum in the $Z$-diagonal block of $Q$ is
$(n-1)+(k-1)=n+k-2$; that in the $A$-diagonal block is
$(n-2)+(a-1)=n+t-k-4$. The other row sums follow directly from
the adjacencies, yielding \eqref{eq:N}.
This quotient is nonnegative and irreducible. Its positive
eigenvector lifts to a positive eigenvector of $Q(H_{n,k,t})$,
so its Perron eigenvalue is $q(H_{n,k,t})$.
It is also similar, by the diagonal matrix of square roots of the
cell sizes, to a real symmetric matrix, so all its eigenvalues are
real.
\end{proof}

\section{Laplacian extremality and complete spectra}\label{sec:laplacian}

We recall the complement identity and its equality consequence,
which are classical \cite[Theorem 3.6 and Corollary 3.7]{Mohar1991}.
They give the entire Laplacian extremal problem here without
requiring a further structural classification.

\begin{theorem}\label{thm:Lmax}
Let $G$ be a graph of order $n\ge2$. Then
$\mu_1(G)\le n$, and
\begin{equation}\label{eq:Lmultiplicity}
 \mult_{L(G)}(n)=c(\overline G)-1,
\end{equation}
where $c(\overline G)$ is the number of components of the complement.
Consequently,
\[
 \mu_1(G)=n\quad\Longleftrightarrow\quad
 \overline G\text{ is disconnected}.
\]
\end{theorem}
\begin{proof}
The degree matrices of $G,\overline G$ sum to $(n-1)I_n$,
and their adjacency matrices sum to $J-I_n$. Hence
\begin{equation}\label{eq:Lcomplement}
 L(G)+L(\overline G)=nI_n-J.
\end{equation}
Both Laplacians are positive semidefinite, since
$x^{\mathsf T}L(G)x=\sum_{uv\in E(G)}(x_u-x_v)^2$.
On $\ones^\perp$, \eqref{eq:Lcomplement} gives
$x^{\mathsf T}L(G)x\le n\|x\|^2$.
On the span of $\ones$, $L(G)$ is zero. These invariant orthogonal
subspaces prove $\mu_1(G)\le n$.

If $L(G)x=nx$, multiplication by $\ones^{\mathsf T}$
gives $\ones^{\mathsf T}x=0$, and \eqref{eq:Lcomplement}
then gives $L(\overline G)x=0$. The converse follows from the
same identity. Therefore the eigenspace for $n$ is
\[
 \ker L(\overline G)\cap\ones^\perp.
\]
The Laplacian quadratic form shows that
$\ker L(\overline G)$ consists exactly of vectors constant on each
component of $\overline G$. It has dimension $c(\overline G)$,
and the condition of being orthogonal to $\ones$ removes one
dimension. This proves \eqref{eq:Lmultiplicity} and the equality
criterion.
\end{proof}

\begin{corollary}\label{cor:Lextremal}
For $2\le t\le n$, the maximum of $\mu_1(G)$ over connected graphs
of order $n$ and integrity $t$ is $n$.
For $1\le k<t-1\le n-2$, the same maximum is $n$ when
$\kappa(G)\le k$ is imposed as well.
In either class a graph is extremal if and only if its complement
is disconnected. The same conclusions hold with $\I(G)\le t$.
\end{corollary}
\begin{proof}
Theorem~\ref{thm:Lmax} supplies the upper bound and equality
criterion. Proposition~\ref{prop:feasibility} supplies the feasible
graphs $F_{n,t}$ and $H_{n,k,t}$. Each has a universal vertex,
so its complement is disconnected, and each attains $n$.
They are feasible also for the upper-bound versions of the classes.
\end{proof}

For a multiset of eigenvalues, $\lambda^{[a]}$ denotes $a$ copies
of $\lambda$; an exponent zero means that the term is absent.

\begin{theorem}\label{thm:Fspectrum}
For $2\le t\le n$,
\begin{equation}\label{eq:Fspectrum}
 \Spec_L(F_{n,t})=
 \bigl\{0,\ n^{[t-1]},\ (t-1)^{[n-t]}\bigr\}.
\end{equation}
\end{theorem}
\begin{proof}
Put $r=t-1$ and $p=n-t+1$. A vector supported on the
$r$-clique with coordinate sum zero is an eigenvector of $L$
with eigenvalue $n$; these vectors have dimension $r-1$.
A vector supported on the independent set with sum zero has
eigenvalue $r$, giving dimension $p-1$.
On the remaining two-dimensional space, consisting of vectors
constant on each part, the quotient is
\[
 \begin{pmatrix}p&-p\\-r&r\end{pmatrix}.
\]
Its eigenvalues are $0$ and $r+p=n$.
These invariant subspaces have dimensions $(r-1)+(p-1)+2=n$.
Thus the total multiplicity of $n$ is $r=t-1$, proving
\eqref{eq:Fspectrum}, including the boundary cases $r=1$ or $p=1$.
\end{proof}

\begin{lemma}[Laplacian spectrum of a join]\label{lem:Ljoin}
Let $X,Y$ have orders $u,v\ge1$, and list their Laplacian
eigenvalues as
$0=\alpha_1,\alpha_2,\ldots,\alpha_u$ and
$0=\beta_1,\beta_2,\ldots,\beta_v$, with multiplicity.
Then
\[
 \Spec_L(X\vee Y)=
 \{0,u+v,\ v+\alpha_2,\ldots,v+\alpha_u,\
                 u+\beta_2,\ldots,u+\beta_v\}.
\]
Only one copy of zero is omitted from each component list, even if
$X$ or $Y$ is disconnected.
\end{lemma}
\begin{proof}
In the natural vertex order,
\[
 L(X\vee Y)=
 \begin{pmatrix}L(X)+vI_u&-J_{u,v}\\
                -J_{v,u}&L(Y)+uI_v\end{pmatrix}.
\]
Choose orthonormal eigenbases of $L(X),L(Y)$ whose first vectors
are the normalized all-ones vectors. Eigenvectors orthogonal to
the all-ones vector in the first part extend by zero and acquire
the eigenvalue shift $v$; those in the second part acquire the
shift $u$. The remaining two-dimensional constant-on-parts space
has quotient
$\left(\begin{smallmatrix}v&-v\\-u&u\end{smallmatrix}\right)$,
with eigenvalues $0,u+v$. The dimensions total $u+v$,
which proves the assertion. This is the standard join formula
\cite[Corollary 3.7]{Mohar1991}.
\end{proof}

\begin{theorem}\label{thm:Hspectrum}
Let $1\le k<t-1\le n-2$. Then
\begin{equation}\label{eq:Hspectrum}
 \Spec_L(H_{n,k,t})=
 \bigl\{0,\ k,\ n^{[k]},\
 (n-1)^{[t-k-1]},\
 (t-1)^{[n-t-1]}\bigr\}.
\end{equation}
In particular, the multiplicity of $n$ is exactly $k$.
\end{theorem}
\begin{proof}
Put $a=t-k-1\ge1$, $p=n-t\ge1$, and
$X=K_a\vee\overline K_p$. By Theorem~\ref{thm:Fspectrum},
\[
 \Spec_L(X)=\{0,(a+p)^{[a]},a^{[p-1]}\}.
\]
Thus $Y=K_1\cup X$ has two zero eigenvalues, together with the
listed nonzero eigenvalues of $X$. Now
$H_{n,k,t}=K_k\vee Y$, where $|Y|=n-k$.
Lemma~\ref{lem:Ljoin} gives one zero, one eigenvalue $n$,
$k-1$ more eigenvalues $n$, one eigenvalue $k$ from the second
zero of $Y$, and the shifted eigenvalues
\[
 k+(a+p)=n-1\quad\text{with multiplicity }a,
 \qquad k+a=t-1\quad\text{with multiplicity }p-1.
\]
This gives \eqref{eq:Hspectrum}. None of its other displayed
eigenvalues equals $n$, so the multiplicity is exactly $k$.
The total multiplicity is
$1+1+k+a+p-1=n$.
\end{proof}

\begin{remark}
In the four-cell order used in Proposition~\ref{prop:qH}, the
Laplacian quotient is
\[
 \begin{pmatrix}
 n-k&-1&-a&-p\\
 -k&k&0&0\\
 -k&0&k+p&-p\\
 -k&0&-a&k+a
 \end{pmatrix}.
\]
Its eigenvalues are $0,k,n-1,n$. Indeed, the within-cell
zero-sum spaces contribute
$n^{[k-1]},(n-1)^{[a-1]},(t-1)^{[p-1]}$; subtracting these
from \eqref{eq:Hspectrum} leaves precisely the four stated
quotient eigenvalues. Thus \eqref{eq:Hspectrum} also evaluates
the quotient explicitly.
\end{remark}

\section{Adjacency formulas and asymptotics}\label{sec:adjacency}

We finish with exact formulas and fixed-parameter expansions for
the two constructions. All implicit constants below may depend
on the fixed integer parameters, but not on the growing order.

\begin{proposition}\label{prop:adjF}
For $2\le t\le n$,
\begin{equation}\label{eq:adjF}
 \rho(F_{n,t})=
 \frac{t-2+\sqrt{(t-2)^2+4(t-1)(n-t+1)}}{2}.
\end{equation}
For fixed $t$, as $n\to\infty$,
\begin{align*}
 \rho(F_{n,t})
 &=\sqrt{(t-1)(n-t+1)}+\frac{t-2}{2}\\
 &\quad+\frac{(t-2)^2}{8\sqrt{(t-1)(n-t+1)}}
       +O_t(n^{-3/2}).
\end{align*}
\end{proposition}
\begin{proof}
The adjacency quotient for the clique and independent set is
$\left(\begin{smallmatrix}t-2&n-t+1\\t-1&0\end{smallmatrix}\right)$.
Its positive eigenvalue is the Perron eigenvalue, giving
\eqref{eq:adjF}. With $r=t-1$ and $P=n-t+1$, expand
\[
 \sqrt{(r-1)^2+4rP}
 =2\sqrt{rP}\left(1+\frac{(r-1)^2}{8rP}
                    +O_r(P^{-2})\right).
\]
Substitution proves the expansion.
\end{proof}

\begin{lemma}[Resolvent for an attached vertex]\label{lem:resolvent}
Let $G_0=K_r\vee\overline K_p$, where $r,p\ge1$.
Let $A_0=A(G_0)$, let $\lambda_0=\rho(G_0)$, and let $u$
be its positive unit Perron vector. Let $v$ be the indicator
vector of $k$ vertices of its clique, where $1\le k\le r$.
Attach a new vertex with neighbourhood given by $v$, producing
a graph $G$ with $\lambda=\rho(G)$. Then $\lambda>\lambda_0$,
and
\begin{equation}\label{eq:resolvent}
 \lambda=\frac{(v^{\mathsf T}u)^2}{\lambda-\lambda_0}
          +R(\lambda),\qquad
 0\le R(\lambda)\le\frac{k}{\lambda}.
\end{equation}
Moreover,
\begin{equation}\label{eq:exact-resolvent}
 \lambda=
 \frac{k-k^2/r}{\lambda+1}
 +\frac{k^2}{r}\,
   \frac{\lambda}{\lambda^2-(r-1)\lambda-rp}.
\end{equation}
\end{lemma}
\begin{proof}
In the old-vertices-first order,
\[
 A(G)=\begin{pmatrix}A_0&v\\v^{\mathsf T}&0\end{pmatrix}.
\]
For $\varepsilon>0$, evaluate its Rayleigh quotient at
$(u,\varepsilon)^{\mathsf T}$. The quotient is
\[
 \frac{\lambda_0+2\varepsilon v^{\mathsf T}u}
      {1+\varepsilon^2}>\lambda_0
\]
for all sufficiently small positive $\varepsilon$, because
$v^{\mathsf T}u>0$. Hence $\lambda>\lambda_0$.
If $(y,z)^{\mathsf T}$ is a positive Perron vector of $A(G)$,
then
$(\lambda I-A_0)y=vz$ and $\lambda z=v^{\mathsf T}y$.
Since $\lambda>\lambda_0$ and $z>0$, elimination gives
\[
 \lambda=v^{\mathsf T}(\lambda I-A_0)^{-1}v.
\]

For completeness, the adjacency spectrum of $G_0$ is
\[
 \{\lambda_0,\lambda_-,
             (-1)^{[r-1]},0^{[p-1]}\},\qquad
 \lambda_\pm=
 \frac{r-1\pm\sqrt{(r-1)^2+4rp}}{2},
\]
where $\lambda_+=\lambda_0$ and $\lambda_-<0$.
The two roots come from the constant-on-parts space; the other
eigenvalues come from the within-part zero-sum spaces.
Expanding the resolvent in an orthonormal eigenbasis and separating
the Perron term gives \eqref{eq:resolvent}.
All other eigenvalues are nonpositive, so the remaining
denominators are at least $\lambda$. Their nonnegative
numerators sum to at most $\|v\|^2=k$, proving the bounds on $R$.

To obtain \eqref{eq:exact-resolvent}, let
$e=r^{-1/2}\ones_{K_r}$ and
$f=p^{-1/2}\ones_{\overline K_p}$, extended by zero off their
respective parts. Decompose
$v=(k/\sqrt r)e+w$, where $w$ is supported on the clique,
orthogonal to $e$, and
$\|w\|^2=k-k^2/r$. On this zero-sum space $A_0=-I$.
On $\operatorname{span}\{e,f\}$, its symmetric matrix is
\[
 \begin{pmatrix}r-1&\sqrt{rp}\\\sqrt{rp}&0\end{pmatrix}.
\]
The upper-left entry of its resolvent is
$\lambda/(\lambda^2-(r-1)\lambda-rp)$.
The orthogonal decomposition of $v$ now proves
\eqref{eq:exact-resolvent}.
\end{proof}

\begin{theorem}\label{thm:Hasymp}
Fix integers $1\le k<t-1$, and put $r=t-1$ and $p=n-t$.
As $n\to\infty$,
\begin{equation}\label{eq:Hasymp}
 \rho(H_{n,k,t})
 =\sqrt{rp}+\frac{r-1}{2}
  +\frac{(r-1)^2/8+k^2/(2r)}{\sqrt{rp}}
  +O_{r,k}(p^{-1}).
\end{equation}
In particular,
\[
 \rho(H_{n,k,t})
 =\sqrt{(t-1)(n-t)}+\frac{t-2}{2}
       +O_{k,t}((n-t)^{-1/2}),
 \qquad
 \lim_{n\to\infty}\frac{\rho(H_{n,k,t})}{\sqrt n}
 =\sqrt{t-1}.
\]
\end{theorem}
\begin{proof}
The graph $H_{n,k,t}$ is obtained from $G_0=K_r\vee\overline K_p$
as in Lemma~\ref{lem:resolvent}. Write
$\lambda=\rho(H_{n,k,t})$ and $\delta=\lambda-\lambda_0>0$.
The exact formula for $\lambda_0$ gives
\begin{equation}\label{eq:lambda0expansion}
 \lambda_0=\sqrt{rp}+\frac{r-1}{2}
       +\frac{(r-1)^2}{8\sqrt{rp}}+O_r(p^{-3/2}).
\end{equation}
By \eqref{eq:resolvent},
\[
 \delta=\frac{(v^{\mathsf T}u)^2}{\lambda-R(\lambda)}.
\]
For sufficiently large $p$, the denominator is at least
$\lambda-k/\lambda>0$. Since
$(v^{\mathsf T}u)^2\le\|v\|^2\|u\|^2=k$ and
$\lambda\ge\lambda_0\asymp_r\sqrt p$, it follows that
\[
 0<\delta\le\frac{k}{\lambda-k/\lambda}
           =O_{r,k}(p^{-1/2}).
\]

Let $\alpha$ be the common coordinate of $u$ on the clique, and
$\beta$ its coordinate on the independent set.
The eigenvector equation gives $\beta=r\alpha/\lambda_0$, while
normalization gives $r\alpha^2+p\beta^2=1$.
Since $\lambda_0^2-(r-1)\lambda_0-rp=0$, these identities imply
\[
 \alpha^2
 =\frac{1}{r(1+rp/\lambda_0^2)}
 =\frac{\lambda_0}{r(2\lambda_0-r+1)}.
\]
Consequently
\[
 (v^{\mathsf T}u)^2=k^2\alpha^2
  =\frac{k^2}{2r}+O_{r,k}(p^{-1/2}).
\]
We already know
$\lambda=\sqrt{rp}+O_{r,k}(1)$ and
$R(\lambda)=O_{r,k}(p^{-1/2})$. Therefore
\[
 \frac{1}{\lambda-R(\lambda)}
  =\frac{1}{\sqrt{rp}}+O_{r,k}(p^{-1}),
\]
and hence
\[
 \delta=\frac{k^2}{2r\sqrt{rp}}+O_{r,k}(p^{-1}).
\]
Adding this to \eqref{eq:lambda0expansion} proves
\eqref{eq:Hasymp}. The two final conclusions follow by retaining
the first two terms and using $p=n-t$ with $t$ fixed.
\end{proof}

\section{Conclusion}\label{sec:conclusion}

We have determined the signless Laplacian spectral extremizers for
connected graphs with bounded vertex integrity, both without an
additional connectivity restriction and with a prescribed upper
bound on vertex-connectivity. The extremal graphs are, respectively,
$F_{n,t}=K_{t-1}\vee\overline K_{n-t+1}$ and
$H_{n,k,t}=K_k\vee(K_1\cup
(K_{t-k-1}\vee\overline K_{n-t}))$. In both cases the integrity is
exactly $t$, which shows that the upper-bound and exact-integrity
formulations have the same extremizer. The bounded-connectivity
argument is based on choosing integrity and cut certificates
simultaneously, making them nested, and then using strict
Perron-vector edge relocations and saturation to force the extremal
structure.

The results complement the adjacency-spectral investigation of
Dai, Shen, Song, and Liu \cite{DaiShenSongLiu2027}. The two spectral
problems are genuinely distinct: the signless Laplacian includes
degree information, and the proof here uses its quadratic form
$\sum_{uv\in E(G)}(x_u+x_v)^2$. The quotient matrices provide exact
implicit descriptions of the signless Laplacian radii, while the
Laplacian spectra and adjacency asymptotics describe further
spectral features of the extremal constructions.

Several related questions remain open. One may seek a direct
comparison between the adjacency and signless Laplacian extremizers
for other integrity-connectivity parameter ranges, investigate
stability versions describing graphs whose signless Laplacian
spectral radius is close to the maximum, or replace the upper bound
on vertex-connectivity by an edge-connectivity constraint. Such
extensions would require arguments that retain control of both
integrity certificates and the relevant connectivity obstruction.

\end{document}